\documentclass{article}
\usepackage[utf8]{inputenc}
\usepackage{comment, amsmath, amssymb, amsthm, authblk}
\usepackage[mathlines]{lineno}

\newtheorem{lem}{Lemma}[section]
\newtheorem{cor}{Corollary}[section]
\newtheorem{pro}{Proposition}[section]
\newtheorem{teo}{Theorem}[section]
\theoremstyle{definition}

\newcommand{\T}[1]{\Tilde{#1}}
\newcommand{\mc}[1]{\mathcal{#1}}

\title{Density of expansivity for geodesic flows of compact higher genus surfaces without conjugate points}
\author[1]{Edhin F. Mamani\thanks{emamani@ufmg.br}}
\author[2]{Rafael Ruggiero\thanks{rorr@mat.puc-rio.br}}
\affil[1]{Instituto de Ciências Exatas ICEx, Universidade Federal de Minas Gerais, Av. Antônio Carlos 6627, Belo Horizonte 31270-901, Brazil.}
\affil[2]{Departamento de Matemática, Pontifícia Universidade Católica do Rio de Janeiro, Rua Marquês de São Vicente 225, Rio de Janeiro 22451-900, Brazil.}
\date{}

\begin{document}
	
	\maketitle
	
	\begin{abstract}
	   Let $(M,g)$ be a compact connected $C^{\infty}$ surface without conjugate points of genus greater than one. We show that set of geodesics without strips forms a dense set of orbits in the unit tangent bundle. This fact was known assuming no focal points as a consequence of a result of Coudène and Shapira. They showed that flat strips are periodic and hence form a set of zero measure in the unit tangent bundle.
	\end{abstract}
	
	\section{Introduction}
 
The study of the geodesic flow of compact surfaces without conjugate points of genus greater than one has had major advances in recent years. Since Morse's pioneering work in the 1920's \cite{morse24}, such geodesic flows are considered paradigms of non-uniformly hyperbolic conservative dynamics. Morse essentially shows that geodesic flows of compact higher genus surfaces without conjugate points can be regarded as coarse Anosov flows. Gromov \cite{grom87} and Ghys \cite{ghys84} independently showed that such flows are topologically semi-conjugate to geodesic flows of hyperbolic surfaces, i.e., surfaces of constant negative curvature. Namely, there exists a continuous surjective map between the correspondent unit tangent bundles sending orbits of the geodesic flow to orbits of the geodesic flow of hyperbolic surfaces. This semi-conjugacy fails to be a conjugacy because there might exist infinite many orbits mapping into single orbits of the hyperbolic geodesic flow. These sets of infinite orbits form strips of bi-asymptotic orbits of the geodesic flow.
	

The outstanding developments of hyperbolic dynamics theory in the 1970's and 1980's raised many natural questions about topological and smooth dynamics of geodesic flows of compact surfaces without conjugate points. Just to mention some of them: Are there non-zero Lyapunov exponents? How ``expansive" is the dynamics? Is the measure of maximal entropy unique?. For topological reasons, it was already knew that topological entropy of geodesic flows of compact higher genus surfaces is always positive. Then, Katok's work \cite{katok80} in the late 1970's for surface diffeomorphisms implies the existence of a hyperbolic invariant measure for the geodesic flow. In particular, there always exist non-zero Lyapunov exponents and invariant measures with positive metric entropy by Pesin's formula, regardless of the no conjugate points assumption. 

Knieper \cite{knie98} in the late 1990's showed that the geodesic flow of a compact rank-1 manifold of non-positive curvature has a unique measure of maximal entropy. This was striking from the perspective of Bowen's work about uniqueness of the measure of maximal entropy for expansive homeomorphisms having the specification property \cite{bowen71endo}. Knieper's approach is based on the theory of Patterson-Sullivan measure. This viewpoint led to the uniqueness of the measure of maximal entropy for geodesic flows of compact higher genus surfaces without conjugate points by Knieper-Climenhaga-War \cite{clim21}. 
 
Concerning the topological dynamics of the geodesic flow, Coudène and Shapira \cite{cou14} showed that every non-trivial strip is periodic assuming nonpositive curvature. Thus, the set of strips is countable for compact higher genus surfaces with non-positive curvature. This in turn yields that the set of ``expansive points", points without strips, is an open set with total measure. Gelfert-Ruggiero \cite{gelf19,gelf20} and Mamani \cite{mam23} showed that geodesic flow of a compact higher genus surface without conjugate points is time-preserving semi-conjugate to a continuous expansive flow acting on a compact metric space (a smooth surface if Green bundles are continuous). This result is somehow surprising from the viewpoint of spectral rigidity theory \cite{croke90,croke92}, since time-preserving semi-conjugations preserve length spectrum of the metric. Moreover, a byproduct of this result gives an alternative proof of the uniqueness of the measure of maximal entropy applying Bowen's theory without using Patterson-Sullivan theory.  
	
The main result of the article deepens the knowledge about how much expansive is the geodesic flow of compact higher genus surfaces without conjugate points. This extends in many senses the consequences of Coudène-Shapira's result about expansiveness: 
		
	\begin{teo} \label{main}
	Let $(M,g)$ be a compact connected $C^{\infty}$ surface without conjugate points of genus greater than one. Then the set of expansive points of the geodesic flow is dense in the unit tangent bundle. 
	\end{teo}

A precise definition of expansive points is given in Section 2. Theorem \ref{main} has many interesting consequences. For example, density of periodic orbits of the geodesic flow (already proved by Knieper-Climenhaga-War \cite{clim21}) and expansivity of the geodesic flow when restricted to a dense set. 

Let us give a brief description of the paper. The main difficulty to prove Theorem \ref{main} is the lack of information about the geometry of strips. Coudène-Shapira's work is based on the flat strip theorem \cite{pesin77} which no longer holds for surfaces without conjugate points and no further assumptions. Though, J-P. Schroeder \cite{roder17} extends this property to strips with positive inner width. 

So instead of dealing with strip geometry we follow the ideas of Gelfert-Ruggiero's work \cite{gelf20,mam23}, where geodesics flows of compact higher genus surfaces without conjugate points are regarded as expansive flows ``up to strips". A strip is generated by the geodesic flow action on a bi-asymptotic class. This class is defined as the intersection between the stable and unstable sets and is a compact connected curve (Subsection 2.3). After identifying these classes with single points we get an expansive flow. 

The key results of the article are in Sections 4 and 5. There we show that return maps of recurrent points have a weak hyperbolic behavior when restricted to special foliated neighborhoods of a bi-asymptotic class. From this we deduce in Section 6 that there cannot exist open sets foliated by bi-asymptotic classes. Moreover, every open set of the unit tangent bundle contains at least one point which agrees with its bi-asymptotic class, i.e., an expansive point.

In Section 7 we give two applications of Theorem \ref{main}: the measure of maximal entropy has full support and periodic points are dense, facts already proved in \cite{clim21} by different methods. We finish the article in Section 8 with a discussion about the relationship between expansive points, expansivity of the geodesic flow and topological transversality of horospheres. We propose a sort of topological version of the well-known Eberlein's characterization of Anosov geodesic flows using transversality of Green bundles.
 	
\section{Preliminaries}
	
\subsection{Compact manifolds without conjugate points}\label{m}
	Let us begin  with basic definitions and notations that we shall use throughout the paper. Let $(M,g)$ be a $C^{\infty}$ compact connected, boundaryless Riemannian manifold, $TM$ be its tangent bundle and $T_1M$ be its unit tangent bundle. Let $\tilde{M}$ be the universal covering of $M$, let $\pi:\tilde{M}\to M$ be the covering map, and let $d\pi:T\tilde{M}\to TM$ the natural covering projection. The universal covering $(\tilde{M},\tilde{g})$ is a complete Riemannian manifold with the pullback metric $\tilde{g}=\pi^*g$. 
	
	The canonical projection is the map $P : TM \longrightarrow M$ given by $P(p,v) = p$. 
	
	A manifold $M$ has no conjugate points if the exponential map $\exp_p$ is non-singular at every $p\in M$. In particular, $\exp_p$ is a covering map for every $p\in M$ (p. 151 of \cite{doca92}).
	
	Denote by $\nabla$ the Levi-Civita connection associated to $g$. A geodesic is a smooth curve $\gamma\subset M$ with $\nabla_{\dot{\gamma}}\dot{\gamma}=0$. For every $\theta=(p,v)\in TM$, denote by $\gamma_{\theta}$ the unique geodesic with initial conditions $\gamma_{\theta}(0)=p$ and $\dot{\gamma}_{\theta}(0)=v$. The geodesic flow $\phi_t$ is defined by 
	\[ \phi: \mathbb{R}\times TM\to TM \qquad (t,\theta)\mapsto \phi_t(\theta)=(\gamma_{\theta}(t), \dot{\gamma}_{\theta}(t)). \]
	All geodesics will be parametrized by arc-length, so we shall restrict the geodesic flow to $T_1M$.
	
	The restriction $P:T_{1}M \longrightarrow M$ at $(p,v) \in T_{1}M$ of the canonical projection to $T_{1}M$ gives rise to the following objects: the vertical fiber $V_{(p,v)}$ defined by 
	$$V(p,v) = P^{-1}(p)= \{ (p,v), \mbox{ }v \in T_{p}M, \mbox{ }\parallel v \parallel = 1\}$$ 
	and the vertical subspace $\mathcal{V}_{(p,v)}$ that is the kernel of the differential of the map $P:T_{1}M \longrightarrow M$ at $(p,v) \in T_{1}M$. 
	
	We shall consider the metric structure in $T_{1}M$ defined by the Sasaki metric $d_{S}$ associated to $g$ (See for details Section 1.3 of (reference)). The Sasaki metric comes from an inner product defined in $TT_{1}M$, and the canonical projection 	$P : (T_{1}M, d_{g}) \longrightarrow (M,g)$ is a Riemannian submersion. We shall often omit the sub-index $g$ in $d_{g}$ to simplify the notation, the ambient space under consideration will uniquely determine the metric.  
	
	\subsection{Busemann functions and horospheres}\label{h}
	
	Let us briefly introduce some important objects in the universal covering which give a good description of the global geometry of geodesics. We follow \cite{esch77} and part II of \cite{pesin77}. Let $\theta\in T_1\tilde{M}$ and $\gamma_{\theta}$ be the geodesic induced by $\theta$. We define the forward Busemann function by
	\[ b_{\theta}: \tilde{M}\to \mathbb{R} \qquad p\mapsto b_{\theta}(p)=\lim_{t\to \infty}d(p,\gamma_{\theta}(t))-t. \]
	Given $\theta=(p,v)\in T_1\Tilde{M}$ let us denote $-\theta:=(p,-v)\in T_1\Tilde{M}$. The stable and unstable horosphere of $\theta$ are defined by
	\[  H^+(\theta)=b_{\theta}^{-1}(0)\subset \tilde{M} \quad \text{ and }\quad  H^-(\theta)=b_{-\theta}^{-1}(0)\subset \tilde{M}. \]
	We lift these horospheres to $T_1\tilde{M}$. Denote by $\nabla b_{\theta}$ the gradient vector field of $b_{\theta}$. We define the stable and unstable horocycle of $\theta$ by
	\[  \mathcal{\tilde{F}}^s(\theta)=\{ (p,-\nabla_pb_{\theta}): p\in H^+(\theta) \}\quad \text{ and }\quad \mathcal{\tilde{F}}^u(\theta)= \{ (p,\nabla_pb_{-\theta}): p\in H^-(\theta) \}. \]
Note that $-\nabla_pb_{\theta}$ forms a vector field on $T_1\T{M}$. The integral curves of $-\nabla_pb_{\theta}$ are called Busemann asymptotes of $\gamma_{\theta}$. The horocycles project onto the horospheres by the canonical projection $\tilde{P}$. For every $\theta\in T_1\tilde{M}$, we define the stable and unstable families of horocycles by
	\[ \mathcal{\tilde{F}}^s=(\mathcal{\tilde{F}}^s(\theta) )_{\theta\in T_1\tilde{M}} \quad \text{ and }\quad \mathcal{\tilde{F}}^u=( \mathcal{\tilde{F}}^u(\theta) )_{\theta\in T_1\tilde{M}}. \]
	We also define the center stable and center unstable sets of $\theta$ by 
	\[\mathcal{\tilde{F}}^{cs}(\theta)=\bigcup_{t\in \mathbb{R}} \mathcal{\tilde{F}}^s(\phi_t(\theta))  \quad \text{ and }\quad \mathcal{\tilde{F}}^{cu}(\theta)=\bigcup_{t\in \mathbb{R}} \mathcal{\tilde{F}}^u(\phi_t(\theta)).\]
	We can define the above objects in the case of $T_1M$. For every $\theta\in T_1M$, 
	\[  \mathcal{F}^*(\theta)=d\pi (\mathcal{\tilde{F}}^*(\tilde{\theta}))\subset T_1M \quad \text{ and }\quad \mathcal{F}^*=d\pi (\mathcal{\tilde{F}}^*), \quad *=s,u,cs,cu; \]
	for any lift $\tilde{\theta}\in T_1\tilde{M}$ of $\theta$. In the case of compact surfaces without conjugate points the collections of horospheres have relevant properties.
	\begin{pro}[\cite{esch77,pesin77}]\label{h1}
		Let $M$ be a compact surface without conjugate points of genus greater than one. Then, for every $\theta\in T_1\tilde{M}$,
		\begin{enumerate}
			\item Busemann functions $b_{\theta}$ are $C^{1,L}$ with $L$-Lipschitz unitary gradient for a uniform constant $L>0$ \cite{knip86}.
			\item Horospheres $H^+(\theta),H^-(\theta)\subset \tilde{M}$ and horocycles $\mathcal{\tilde{F}}^s(\theta),\mathcal{\tilde{F}}^u(\theta)\subset T_1\tilde{M}$ are embedded curves.
			\item The families $\mathcal{\tilde{F}}^s,\mathcal{\tilde{F}}^u$ and $\mathcal{F}^s,\mathcal{F}^u$ are continuous foliations of $T_1\tilde{M},T_1M$ respectively, that are invariant by the geodesic flow: for every $t\in \mathbb{R}$,    
			\begin{equation}
				\Tilde{\phi}_t(\mathcal{\tilde{F}}^s(\theta))=\mathcal{\tilde{F}}^s(\Tilde{\phi}_t(\theta)).
			\end{equation}
		\end{enumerate}
	\end{pro}
	
	\subsection{Morse's shadowing and strips}\label{special}
	The celebrated work of Morse \cite{morse24} shows that the global geometry of geodesics in the universal covering of a compact surface without conjugate points and genus greater than one resembles the global geometry of geodesics in the hyperbolic plane. 
	\begin{teo}[\cite{morse24}]\label{mor}
		Let $(M,g)$ be a compact surface without conjugate points and $\tilde{M}$ be its universal covering. Then, there exists $R>0$ such that for every geodesic $\gamma\subset \tilde{M}$ there exists a hyperbolic geodesic $\gamma'\subset \tilde{M}$ with  Hausdorff distance between $\gamma$ and $\gamma'$ bounded above by $R$.
	\end{teo}
	Given two geodesics $\gamma_1,\gamma_2\subset \Tilde{M}$, we say that $\gamma_1$ and $\gamma_2$ are asymptotic if $d(\gamma_1(t),\gamma_2(t))\leq C$ for every $t\geq 0$ and for some $C>0$. If the last inequality holds for every $t\in \mathbb{R}$, $\gamma_1$ and $\gamma_2$ are called bi-asymptotic. Theorem \ref{mor} actually provides a uniform bound for the Hausdorff distance between bi-asymptotic geodesics in $(\tilde{M}, \tilde{g})$. 
	\begin{teo}\label{morse}
		Let $(M,g)$ be a compact surface without conjugate points of genus greater than one. Then there exists $Q(M)>0$ such that the Hausdorff distance between any two bi-asymptotic geodesics is bounded above by $Q(M)$.
	\end{teo}
	For each $\theta\in T_1\tilde{M}$, we define the intersections 
	\[ I(\theta)=H^+(\theta)\cap H^-(\theta)\subset \tilde{M} \quad \text{ and }\quad  \mathcal{\tilde{I}}(\theta)=\mathcal{\tilde{F}}^s(\theta)\cap \mathcal{\tilde{F}}^u(\theta)\subset T_1\tilde{M}.  \]
	We shall call $\mathcal{\tilde{I}}(\theta)$ the bi-asymptotic class of $\theta$, or simply the class of $\theta$. Furthermore, $\tilde{P}(\mathcal{\tilde{I}})=I(\theta)$ where $\tilde{P}:T_1\tilde{M}\to M$ is the canonical projection.
	
	We observe that for every $\eta=(q,w)\in \mathcal{\Tilde{I}}(\theta)$ with $q\in I(\theta)$, the geodesic $\gamma_{\eta}$ is bi-asymptotic to $\gamma_{\theta}$. Let us state some of the properties of $I(\theta)$ and $\mathcal{\tilde{I}}(\theta)$ that will be important in the sequel. 
	\begin{pro}\label{tops}
		Let $M$ be a compact surface without conjugate points of genus greater than one and $\tilde{M}$ be its universal covering. Then, for every $\theta\in T_1\tilde{M}$
		\begin{enumerate}
			\item $I(\theta)$ and $\mathcal{\tilde{I}}(\theta)$ are compact connected curves of $\tilde{M}$ and $T_1\tilde{M}$ respectively (Corollary 3.3 of \cite{riff18}).
			\item $Diam(I(\theta))\leq Q$ and $Diam(\mathcal{\tilde{I}}(\theta))\leq \tilde{Q}$ for some $Q(M),\tilde{Q}(M)>0$.
			\item Two bi-asymptotic classes are either disjoint or coincide. 
		\end{enumerate}
	\end{pro}

    \begin{proof}
    Item 2 is straightforward from Morse's Theorem and item 3 holds because bi-asymptoticity is an equivalence relation.
    \end{proof}

    The expansive points form the so-called expansive set 
	\[   \mathcal{R}_0=\{ \theta \in T_1M: \mathcal{F}^s(\theta)\cap \mathcal{F}^u(\theta)= \{ \theta \} \}. \]
    Its complement is called the non-expansive set and denoted by $\mathcal{R}_0^c$. Moreover, Proposition \ref{tops}(1) says that $\mathcal{I}(\theta)$ is a compact connected set of dimension 1 hence homeomorphic to a non-trivial real interval. Thus, any non-trivial class $\mathcal{I}(\theta)$ has two boundary points included in $\mathcal{F}^s(\theta)$ (or $\mathcal{F}^u(\theta)$).
 
    In the following result, proved in \cite{mam23}, we will consider the subspace topology of $\mc{\T{F}}^s(\theta)$ and $\mc{T{F}}^u(\theta)$ for the conclusions.

\begin{lem}\label{que2}
	Let $(M,g)$ be a compact surface without conjugate points of genus greater than one and $\theta\in T_1\Tilde{M}$ be a non-expansive point. Then the boundary points of $\mathcal{\tilde{I}}(\theta)$ are accumulated by expansive points included in $\tilde{\mathcal{F}}^s(\theta)$. In particular, expansive points are dense in $\mc{\T{F}}^s(\theta)\setminus \bigcup_{\eta\in \mc{R}_0^c}\mc{\T{I}}(\eta)$. A similar statement holds for $\tilde{\mathcal{F}}^u(\theta)$.
\end{lem}
	
\subsection{Visibility manifolds}
    Let $M$ be a simply connected Riemannian manifold without conjugate points. For every $x,y\in M$, denote by $[x,y]$ the geodesic segment joining $x$ to $y$. For $z\in M$ we also denote by $\sphericalangle_z(x,y)$ the angle at $z$ formed by $[z,x]$ and $[z,y]$. We say that $M$ is a visibility manifold if for every $z\in M$ and every $\epsilon>0$ there exists $R(\epsilon, z)>0$ such that 
	\[ \text{ if }x,y\in M \text{ with }d(z,[x,y])>R(\epsilon, z) \quad \text{ then }\quad \sphericalangle_z(x,y)<\epsilon.\]
	If $R(\epsilon,z)$ does not depend on $z$ then $M$ is called a uniform visibility manifold.
	\begin{teo}[\cite{eber72}] \label{Eberlein72}
		If $M$ is a compact surface without conjugate points of genus greater than one then its universal covering is a uniform visibility 2-manifold.
	\end{teo}
	In the 1970's, Eberlein observed that the visibility condition was crucial to extend the theory about the interplay between the geodesic flow dynamics and the action of the fundamental group on the ideal boundary of a hyperbolic manifold to the context of compact manifolds without conjugate points. 
 
	\begin{teo}[\cite{eber72,eber73neg2}]\label{v1}
		Let $M$ be a compact surface without conjugate points of genus greater than one. Then
		\begin{enumerate}
			\item The horospherical foliations $\mathcal{F}^s$ and $\mathcal{F}^u$ are minimal, i.e., each leaf is dense.
			\item The geodesic flow $\phi_t$ is transitive and topologically mixing.
			\item For every $\theta,\xi \in T_1\tilde{M}$ with $\theta\not\in \mathcal{\tilde{F}}^{cu}(\xi)$ there exists $\eta_1,\eta_2\in T_1\tilde{M}$ such that
			\[ \mathcal{\tilde{F}}^s(\theta)\cap \mathcal{\tilde{F}}^{cu}(\xi)=\mathcal{\Tilde{I}}(\eta_1) \quad \text{ and }\quad  \mathcal{\tilde{F}}^s(\xi)\cap \mathcal{\tilde{F}}^{cu}(\theta)=\mathcal{\Tilde{I}}(\eta_2).\]
			\item For every $\theta \in T_1M$ and $\epsilon>0$, there exist $\theta'\in \mathcal{I}(\theta)$ and a periodic point $\xi\in T_1M$ such that $d_s(\theta',\xi)<\epsilon$.        
				\end{enumerate}
	\end{teo}
 
    Item (3) describes a well-known phenomenon in hyperbolic dynamics, heteroclinic relations, which can be also displayed as follows: there exist $t_1,t_2\in \mathbb{R}$ such that
	\[ \mathcal{\tilde{F}}^{cs}(\theta)\cap \mathcal{\tilde{F}}^u(\xi)=\mathcal{\Tilde{I}}(\Tilde{\phi}_{t_1}(\eta_1)) \quad \text{ and }\quad \mathcal{\tilde{F}}^{cs}(\xi)\cap\mathcal{\tilde{F}}^u(\theta)=\mathcal{\Tilde{I}}(\Tilde{\phi}_{t_2}(\eta_2)).  \]
	Item (4) is equivalent to density of periodic points of the geodesic flow up to bi-asymptotic classes. Actually, Climegnaga-Knieper-War showed that periodic points are dense \cite{clim21}, but for our purposes item (4) will be enough.  
		
 It is also important to highlight the global geometric properties of visibility manifolds. Let $(M,g)$ be a compact manifold without conjugate points and $\tilde{M}$ be its universal covering. Geodesics rays diverge in $\tilde{M}$ if for every $p\in \tilde{M}$, every $\epsilon,A>0$, there exists $T(p,\epsilon,A)$ such that for every geodesics $\gamma,\beta \subset \tilde{M}$ with same base point $p$ and $\angle(\gamma'(0),\beta'(0))\geq \epsilon$, then $d(\gamma(t),\beta(t))\geq A$ for every $t\geq T(p,\epsilon,A)$. We say that geodesic rays diverge uniformly if $T(p,\epsilon,A)$ does not depend on $p$. We say that $M$ is quasi-convex if there exist constants $A, B>0$ such that for every two geodesics 
\[ \gamma : [t_1, t_2] \to \tilde{M}, \qquad   \beta : [s_1,s_2] \to \tilde{M}, \]
the Hausdorff distance satisfies
\[ d_{H}(\gamma , \beta ) \leq A \sup \{ d(\gamma(t_1), \beta(s_1)), d(\gamma(t_2), \beta(s_2))\}+ B . \]
Eberlein already proved the quasi-convexity of a uniform visibility manifold \cite{eber72} and the central stable global behavior of orbits starting at $\mathcal{\tilde{F}}^s(\theta)$. 

 By Morse's work, these features were already known for universal coverings of compact surfaces without conjugate points of genus greater than one. We state the stable global behavior in our setting.
	\begin{pro}[\cite{eber72}]\label{quasi}
		Let $M$ be a compact surface without conjugate points of genus greater than one and $\tilde{M}$ be its universal covering. Then, there exist $A,B>0$ such that for every $\theta\in T_1\tilde{M}$ and every $\eta\in \mathcal{\tilde{F}}^s(\theta)$, 
		\[   d_s(\Tilde{\phi}_t(\theta),\Tilde{\phi}_t(\eta))\leq Ad_s(\theta,\eta)+B, \quad \text{ for every } t\geq 0.  \]
	\end{pro}

We observe that visibility manifolds belong to a bigger class of manifolds which have the following important property.
\begin{teo}[Lemma 2.9 of \cite{riff18}]\label{asim}
Let $M$ be a compact manifold without conjugate points and $\tilde{M}$ be its universal covering. If $\tilde{M}$ is quasi-convex and geodesic rays diverge uniformly then for every $\theta\in T_1\tilde{M}$, a geodesic $\beta$ is asymptotic to $\gamma_{\theta}$ if and only if $\beta$ is a Busemann asymptote of $\gamma_{\theta}$.
\end{teo}

Finally, we see that an analogous property to Proposition \ref{tops} holds for visibility manifolds.
\begin{teo} \label{connected-class} (\cite{riff18})
Let $(M,g)$ be a compact manifold without conjugate points and $\T{M}$ be its universal covering. If $\tilde{M}$ is quasi-convex and geodesic rays diverge then $\mathcal{\T{I}}(\tilde{\theta})$ is a connected set for every $\tilde{\theta} \in T_{1}\tilde{M}$. 
\end{teo}
	
\section{Foliated neighborhoods and cross sections for the geodesic flow}\label{fol}

The purpose of the section is to define certain foliated open sets and cross sections for the geodesic flow. These objects were introduced by Gelfert and Ruggiero \cite{gelf19}, and showed to be very convenient to study the geodesic flow dynamics. Let us start with some notations. 
 
For every $\tilde{\theta}=(p,v) \in T_{1} \tilde{M}$, recall that $V(\T{\theta}) = \{ (p,w) \in T_{p}\T{M}: \|w\| = 1\}$ and denote by $V_{\delta}(\tilde{\theta}) \subset V(\tilde{\theta})$ a small open ball of radius $\delta>0$ in $V(\tilde{\theta})$ centered at $\tilde{\theta}$, with respect to Sasaki's metric. For $\tau>0$, define 
\[  Q(\tilde{\theta},\delta,\tau) = \bigcup_{|t|<\tau}\T{\phi}_{t}(V_{\delta}(\tilde{\theta})). \]
Let $U^{s}(\tilde{\theta})$ be a connected, open subset of $\tilde{\mathcal{F}}^{s}(\tilde{\theta})$, containing $\mc{\T{I}}(\T{\theta})$ and with compact closure. From Proposition \ref{h1}, we recall that collections $\mathcal{\tilde{F}}^s$ and $\mathcal{\tilde{F}}^u$ form continuous foliations of $T_1\T{M}$ invariant by the geodesic flow. The following results are well-known for surfaces, complete proofs for compact $n$-manifolds without conjugate points are found in \cite{pot23}. 

\begin{lem} \label{vertical-neigh}
	Let $(M,g)$ be a compact surface without conjugate points of genus greater than one. Let $U^{s}(\tilde{\theta})$ be a connected, open subset of $\tilde{\mathcal{F}}^{s}(\tilde{\theta})$, containing $\mc{\T{I}}(\T{\theta})$ and with compact closure. Then the set  
	\begin{eqnarray*} 
		Q^s(\tilde{\theta},\delta,\tau,U^{s}(\tilde{\theta})) & = & \bigcup_{\tilde{\xi} \in U^s(\tilde{\theta})} Q(\tilde{\xi},\delta,\tau) \\
		& = & \bigcup_{|t|<\tau, \tilde{\xi}\in U^s(\tilde{\theta})}\T{\phi}_{t}(V_{\delta}(\tilde{\xi}))
	\end{eqnarray*}
	is an open neighborhood of $\tilde{\theta}$ in $T_{1}\tilde{M}$. 
\end{lem}

	\begin{lem} \label{Fol-neigh}
	Let $(M,g)$ be a compact surface without conjugate points of genus greater than one. Then, for every $\tilde{\theta}\in T_{1}\tilde{M}$, $U^{s}(\tilde{\theta})$ as above, there exists $0<\delta'\leq\delta$ such that if $d_s(\tilde{\theta},\tilde{\eta})<\delta'$ we have
	\begin{enumerate}
		\item Each submanifold $\tilde{\mathcal{F}}^{s}(\tilde{\eta})$ crosses each set $Q(\tilde{\xi},\delta,\tau)$ for every 
		$\tilde{\xi}\in U^{s}(\tilde{\theta})$ at just one point $\tilde{\eta}^{s}(\tilde{\xi})$. 
		\item The sets 
    \[  \tilde{\mathcal{F}}^{s}(\tilde{\eta},U^s(\tilde{\theta})) = \tilde{\mathcal{F}}^s(\tilde{\eta})\cap Q^s(\tilde{\theta},\delta,\tau,U^s(\tilde{\theta})) \]
		are all homeomorphic to $U^{s}(\tilde{\theta})$. 
		\item The sets 
		\[  S^{s}(U^s(\tilde{\theta}),\delta') = \bigcup_{\T{\eta} \in V_{\delta'}(\tilde{\theta})}\tilde{\mathcal{F}}^s(\tilde{\eta},U^s(\tilde{\theta}))  \]
		are s-foliated, continuous cross sections for the geodesic flow homeomorphic to $U^s(\tilde{\theta}) \times V_{\delta'}(\tilde{\theta})$, and the sets 
		\[ \Gamma^s(U^s(\tilde{\theta}),\delta',\tau) = \bigcup_{|t|<\tau} \T{\phi}_{t}( S^{s}(U^s(\tilde{\theta}),\delta'))   \]
		are s-foliated open neighborhoods of $\tilde{\theta}$ in $T_1\T{M}$.
	\end{enumerate}
	Similar statements hold to u-foliated sections and open neighborhoods. 
\end{lem}
Similar objects can be defined in $T_{1}M$ since given any $\tilde{\theta} \in T_{1}\tilde{M}$ and $U^s(\tilde{\theta})$ as above, we can choose $\delta,\tau>0$ small enough such that the following sets are homeomorphic to their images under the projection map $d\pi: T_1\tilde{M}\to T_1M$ with $\theta=d\pi(\T{\theta})$.
\[ Q(\theta,\delta,\tau) = d\pi(Q(\tilde{\theta},\delta,\tau) )   \]
\[ U^s(\theta) = d\pi(U^s(\tilde{\theta}))  \]
\[ Q^s(\theta,\delta,\tau,U^s(\theta)) = d\pi ( Q^s(\tilde{\theta}),\delta,\tau,U^s(\tilde{\theta}))  \]
\[  \mathcal{F}^{s}(\eta,U^s(\theta)) = d\pi( \tilde{\mathcal{F}}^s(\tilde{\eta},U^s(\tilde{\theta})) \]
\[  S^s(U^s(\theta),\delta') = d\pi(S(U^s(\tilde{\theta})\delta'))  \]
\[  \Gamma^s(U^s(\theta),\delta',\tau) = d\pi( \Gamma^s(U^s(\tilde{\theta}),\delta',\tau)  .\]

\section{Contracting dynamics in $\mathcal{F}^{s}(\theta)$ out of the bi-asymptotic class $\mathcal{I}(\theta)$}\label{scont}

The goal of the section is to show that dynamics of the geodesic flow restricted to $\mathcal{F}^{s}(\theta)$ behaves as a topological contraction outside of the bi-asymptotic class $\mathcal{I}(\theta)$. The main result of the section is inspired by Theorem 5.2 in Potrie-Ruggiero's paper \cite{pot23}. In this work, a similar statement is proved for bi-asymptotic classes in a neighborhood of a hyperbolic closed geodesic for a compact $n$-manifold without conjugate points assuming divergence of geodesic rays and quasi-convexity. Let us first introduce some notations. 

For every $\T{\theta}\in T_1\T{M}$, we define a homeomorphism that assigns 'coordinates' to points in the s-foliated cross section as follows:
\begin{align*}
    f_{\T{\theta}}: S^s(U^s(\T{\theta}),\delta'))&\longrightarrow U^s(\T{\theta}) \times  V_{\delta'}(\T{\theta})    \\
    \omega &\mapsto (\T{\Pi}^h(\omega),\T{\Pi}^v(\omega)),
\end{align*}
so that $\omega \in Q(h(\omega),\delta,\tau)\cap \mathcal{\T{F}}^s(v(\omega),U^s(\T{\theta}))$. In other words, $\T{\Pi}^h(\omega)$ and $\T{\Pi}^v(\omega)$ are the transversal intersections: $Q(\omega,\delta,\tau)\cap U^s(\T{\theta})$ and $V_{\delta'}(\T{\theta})\cap \mathcal{\T{F}}^s(\omega)$ respectively. So, $\T{\Pi}^h(\omega)$ and $\T{\Pi}^v(\omega)$ are in some sense coordinates of $\omega$ in the product $U^s(\T{\theta}) \times  V_{\delta'}(\T{\theta})$.

In particular we observe that restriction of $\tilde{\Pi}^{h}$ to $\T{\mathcal{F}}^s(\T{\eta},U^s(\T{\theta}))$ gives the homeomorphism between $\T{\mathcal{F}}^s(\T{\eta},U^s(\T{\theta})$ and $U^s(\T{\theta})\subset \tilde{\mathcal{F}}^{s}(\T{\theta})$ considered in Lemma \ref{Fol-neigh}. 

We have an analogous situation in an open neighborhood of $\theta=d\pi(\T{\theta})\in T_1M$, because $d\pi$ is a covering map, as we observed at the end of previous section. So, let 
\[  \Pi^h: S^s(U^s(\theta),\delta')\longrightarrow U^s(\theta) \]
be the projection defined by $\Pi^{h} = d\pi \circ \tilde{\Pi}^{h}$.

Let us now consider a sequence of maps on $S^s(U^s(\theta),\delta')$. Let $\theta\in T_1M$ be a recurrent point, $S^s(U^s(\theta),\delta')$, $U^s(\theta)$ as above and $S^s(U^s(\theta),\delta')$ as the end of previous section. Observe that we can choose a sequence of times
$(t_n)\subset \mathbb{R}$ such that $\phi_{t_n}(\theta)\in S^s(U^s(\theta),\delta')$ and $\phi_{t_n}(\theta)\to \theta$. For every $n\geq 1$, the following Poincaré-type map
\[  \mc{P}_n: D_n\subset S^s(U^s(\theta),\delta')\longrightarrow S^s(U^s(\theta),\delta') \]
is defined for every $\xi\in D_n$ as follows. Using the geodesic flow, project the point $\phi_{t_n}(\xi)$ into $S^s(U^s(\theta),\delta')$. By continuity of the geodesic flow, we can reduce $D_n$ if necessary so that $\mc{P}_n$ is well-defined. Clearly, for large enough natural numbers $n\geq 1$, $\phi_{t_n}(\theta)$ is very close to $\theta$ hence $D_n$ contains $\theta$ and is sufficiently big to our purposes.

The main result of the section is the following:

\begin{pro} \label{contraction}
Let $(M,g)$ be a compact surface without conjugate points of genus greater than one. Let $\theta \in T_{1}M$ be a recurrent point, $U^s(\theta)$ be an open relative neighborhood of $\mathcal{I}(\theta)$ in $\mathcal{F}^{s}(\theta)$ with compact closure, and let $\mc{P}_N$ be the map defined  above for large enough $N\geq 1$. Then, for every compact set $K \subset \mathcal{F}^{s}(\theta)$ with $U^s(\theta)\subset K$, there exists $n_{\theta}\geq N$ such that $\Pi^{h}\circ \mathcal{P}_n(K) \subset U^s(\theta)$ for every $n\geq n_{\theta}$. 
\end{pro}

\begin{proof}
By contradiction, suppose there exist a compact set $K_0\subset \mathcal{F}^{s}(\theta)$ containing $U^s(\theta)$ and a sequence $n_k\to \infty$ such that $\Pi^h \circ \mathcal{P}_{n_k}(K_0) \not\subset U^s(\theta)$. Hence, there exists a sequence $\xi_k\in K_0\subset \mathcal{F}^s(\theta)$ such that 
\begin{equation}\label{not}
    \Pi^h \circ \mathcal{P}_{n_k}(\xi_k) \not\in U^s(\theta).
\end{equation}
Now, let $\T{\theta}, \T{\xi}_k\in T_1\T{M}$ be lifts of $\theta,\xi_k$ such that $\T{\xi}_k\in \mc{\T{F}}^s(\T{\theta})$ and $\T{\xi}_k,\T{\theta}\in \T{K}_0$ with $d\pi(\T{K}_0)=K_0$ for some compact set $\T{K}_0\subset T_1\T{M}$. For every $k\geq 1$, we consider $\T{\eta}_k=\phi_{t_{n_k}}(\T{\xi}_k)$ and $\T{\theta}_k=\phi_{t_{n_k}}(\T{\theta})$. Using covering isometries, there exists a subsequence $\T{\eta}_k$ (denoted by same index) converging to some $\T{\eta}\in T_1\T{M}$. Since $\T{\xi}_k\in \mc{\T{F}}^s(\T{\theta})$, we see that 
\begin{equation}\label{con}
    \T{\eta}_k=\phi_{t_{n_k}}(\T{\xi}_k)\in \phi_{t_{n_k}}(\mc{\T{F}}^s(\T{\theta})) =\mc{\T{F}}^s(\phi_{t_{n_k}}(\T{\theta}))=\mc{\T{F}}^s(\T{\theta}_k).
\end{equation}
Thus, by continuity of the foliation $\mc{\T{F}}^s$ we conclude that $\T{\eta}\in \mc{\T{F}}^s(\T{\theta})$. From Equation \eqref{con}, Proposition \ref{quasi} and compactness of $\T{K}_0$ we have
\begin{equation}
    d_s(\phi_t\T{\eta}_k,\phi_t\T{\theta}_k)=d_s(\phi_t\circ \phi_{t_{n_k}}\T{\xi}_k,\phi_t\circ \phi_{t_{n_k}}\T{\theta})\leq Ad_s(\T{\xi}_k,\T{\theta})+B\leq C, 
\end{equation}
for every $t\geq -t_{n_k}$ and some $C>0$. Since $t_{n_k}\to \infty$, this equation implies that $d_s(\phi_t\T{\eta},\phi_t\T{\theta})\leq C$ for every $t\in \mathbb{R}$ hence $\T{\eta}\in \mathcal{\T{I}}(\T{\theta})$ and $\eta=d\pi(\T{\eta})\in \mc{I}(\theta)$. On the other hand, since $\mc{P}_{n_k}(\xi_k)=\phi_{t_{n_k}}(\xi_k)=\eta_k$, Equation \eqref{not} says that $\Pi^h(\eta_k)\not\in U^s(\theta)$. From this, noting that $\theta_k\to \theta$ and $\eta_k\in \mc{F}^s(\theta_k)$, we conclude that $\eta\not\in U^s(\theta)$ and hence $\eta\not\in \mathcal{I}(\theta)$, a contradiction which concludes the proof.

\end{proof}
This proposition yields that the action of the geodesic flow restricted to $\mathcal{F}^{s}(\theta)$ outside the set $\mathcal{I}(\theta)$ is a contraction, the forward orbits of points in every subset of $\mathcal{F}^{s}(\theta)$ containing $\mathcal{I}(\theta)$ must approach the orbits of $\mathcal{I}(\theta)$ whenever $\theta$ is a recurrent point. If $\theta$ is a periodic point, the orbits of the points in $\mathcal{I}(\theta)$ foliate an annulus by closed orbits with the same homotopy class and period. Proposition \ref{contraction} implies that outside $\mathcal{I}(\theta)$ the forward orbits of points in $\mathcal{F}^{s}(\theta)$ approach this annulus.  

\section{The action of the Poincaré-type maps on weak local product neighborhoods}
Let us recall from Lemma \ref{Fol-neigh}, the collection of $s$-foliated cross sections 
\[  S^{s}(U^s(\tilde{\theta}),\delta) = \bigcup_{\T{\eta} \in V_{\delta}(\tilde{\theta})}\tilde{\mathcal{F}}^s(\tilde{\eta},U^s(\tilde{\theta}))  \]
and the collection of $s$-foliated open neighborhoods
\[ \Gamma^s(U^s(\tilde{\theta}),\delta,\tau) = \bigcup_{|t|<\tau} \T{\phi}_{t}( S^{s}(U^s(\tilde{\theta}),\delta))   \]
Moreover, this lemma also says that an analogous construction can be done for the unstable case: there are $u$-foliated cross sections and open neighborhoods containing $\mc{\T{I}}(\T{\theta})$. Namely, let $U^u(\tilde{\theta}) \subset \tilde{\mathcal{F}}^{u}(\tilde{\theta})$ be a relative open neighborhood of $\tilde{\theta}$, $\tilde{\mathcal{F}}^u(\tilde{\eta},U^u(\tilde{\theta}))\subset \tilde{\mathcal{F}}^{u}(\tilde{\eta})$ be the relative open set homeomorphic to $U^u(\tilde{\theta})$ according to an analogous homeomorphism defined in Lemma \ref{Fol-neigh}. So, we denote by
\[  S^{u}(U^u(\tilde{\theta}),\delta) = \bigcup_{\T{\eta} \in V_{\delta}(\tilde{\theta})}\tilde{\mathcal{F}}^u(\tilde{\eta},U^u(\tilde{\theta}))  \]
the u-foliated cross section for the geodesic flow and by
\[ \Gamma^u(U^u(\tilde{\theta}),\delta,\tau) = \bigcup_{|t|<\tau} \T{\phi}_{t}( S^{u}(U^u(\tilde{\theta}),\delta))   \]
the u-foliated open neighborhoods of $\mathcal{\T{I}}(\T{\theta})$. We highlight that $S^{u}(U^u(\tilde{\theta}),\delta)$ and $\Gamma^u(U^u(\tilde{\theta}),\delta,\tau)$ are foliated by open relative sets of horospherical leaves of $\mathcal{\T{F}}^u$. 

We clearly see that $\Gamma^s(U^s(\tilde{\theta}),\delta,\tau)\cap \Gamma^u(U^u(\tilde{\theta}),\delta,\tau)$ is an open neighborhood of $\mc{\T{I}}(\T{\theta})$. So, varying the parameters $U^s(\T{\theta})$, $U^u(\T{\theta})$, $\delta,\tau>0$, we obtain a family of neighborhoods of $\mc{\T{I}}(\T{\theta})$ satisfying 
\[  \bigcap_{U^s(\T{\theta}), U^u(\T{\theta}),\delta,\tau} \Gamma^s(U^s(\tilde{\theta}),\delta,\tau)\cap \Gamma^u(U^u(\tilde{\theta}),\delta,\tau) = \mathcal{\T{I}}(\tilde{\theta}).  \]

Now, for every $\T{\eta}\in S^{s}(U^s(\tilde{\theta}),\delta)$ we define an unstable-type set
\[ \mc{W}^u(\T{\eta})= S^{s}(U^s(\tilde{\theta}),\delta) \cap \mc{\T{F}}^{cu}(\T{\eta}) \subset S^{s}(U^s(\tilde{\theta}),\delta).\]
Furthermore, for every $\delta'<\delta$ we define the set
\[ W(\T{\theta},\delta') = \bigcup_{\T{\eta}\in V_{\delta'}(\T{\theta})} \mc{W}^u(\T{\eta})\subset S^{s}(U^s(\tilde{\theta}),\delta). \]
As a consequence of Theorem \ref{Eberlein72}(3) we obtain the following weak local product.
\begin{lem} \label{crossed-int}
Let $(M,g)$ be a compact surface without conjugate points of genus greater than one, $\T{\theta} \in T_1\T{M}$ and $U^s(\T{\theta})\subset\tilde{\mathcal{F}}^{s}(\tilde{\theta})$ as above. Then, 
	\begin{enumerate}
		\item There exists $\delta>0$ such that for every $\T{\eta},\T{\eta}'\in S^{s}(U^s(\tilde{\theta}),\delta)$ with $d_s(\T{\eta},\T{\eta}')<\delta$, 
  		\[ \tilde{\mathcal{F}}^{s}(\T{\eta}) \cap \tilde{\mathcal{F}}^{cu}(\tilde{\eta}')= \mc{\T{I}}(\T{\xi}), \text{ for some } \T{\xi}\in S^{s}(U^s(\tilde{\theta}),\delta)\]
		\item For every $\delta'< \delta$, the above defined set $W(\T{\theta},\delta')\subset S^{s}(U^s(\tilde{\theta}),\delta)$ is a relative open neighborhood of $\mc{\T{I}}(\T{\theta})$ in $S^{s}(U^s(\tilde{\theta}),\delta)$.
        \item The family $\{ \mc{W}^u(\T{\eta}): \T{\eta}\in V_{\delta'}(\T{\theta}) \}$ is a continuous foliation of $W(\T{\theta},\delta')$ by continuous, $n$-dimensional leaves such that for every $\T{\eta}\in W(\T{\theta},\delta')$ 
		\[  \tilde{\mathcal{F}}^{s}(\tilde{\eta}) \cap \mathcal{W}^{u}(\tilde{\eta})= \mathcal{\T{I}}(\tilde{\eta}).\]
	\end{enumerate} 
	The same statements hold interchanging the indices $s$ and $u$ in the above items.
\end{lem}
\begin{proof}
We sketch the proof. Item (1) is just the heteroclinic relations between orbits given by the visibility condition stated in Theorem \ref{Eberlein72}. Items (2) and (3) are consequences of item (1) and the continuity of horospherical foliations $\mathcal{\T{F}}^s$ and $\mathcal{\T{F}}^u$ given in Lemma \ref{Fol-neigh}. 
\end{proof} 
Let us next consider the projections $\mathcal{W}^{u}(\eta),W(\theta,\delta')\subset S^{s}(U^s(\theta),\delta)$ of $\mathcal{W}^{u}(\tilde{\eta})$, $W(\T{\theta},\delta')\subset S^{s}(U^s(\tilde{\theta}),\delta)$ by the covering map $d\pi$. We can carry the conclusions of Lemma \ref{crossed-int} to this setting. For every $\delta'< \delta$, $W(\theta,\delta')\subset S^{s}(U^s(\theta),\delta)$ is a relative open neighborhood of $\mc{I}(\theta)$ in $S^{s}(U^s(\theta),\delta)$ and the family $\{ \mc{W}^u(\eta): \eta\in V_{\delta'}(\theta) \}$ is a continuous foliation of $W(\theta,\delta')$ by continuous, $n$-dimensional leaves such that for every $\eta\in W(\theta,\delta')$ we have 
		\[  \mathcal{F}_c^{s}(\eta) \cap \mathcal{W}^{u}(\eta)= \mathcal{I}(\eta),\]
where $\mathcal{F}_c^{s}(\eta)$ is the connected component of $\mathcal{F}_c^{s}(\eta)\cap S^{s}(U^s(\theta),\delta)$ containing $\eta$. Notice also that leaves $\mathcal{W}^{u}(\eta)$ might not be subsets of $\mathcal{F}^{u}(\eta)$. 

We next define some projections along unstable leaves as follows
\begin{align*}
    H_{u}: S^{s}(U^s(\theta),\delta) &\longrightarrow V_{\delta}(\theta) \\
    \xi &\mapsto \mathcal{W}^{u}(\xi)\cap V_{\delta}(\theta).
\end{align*}
That is, $H_u(\xi)$ is the projection of $\xi \in S^{s}(U^s(\theta),\delta)$ into $V_{\delta}(\theta)$ along the unstable leaf $\mathcal{W}^{u}(\xi)$. The following result is based on Lemma 5.2 of \cite{pot23}.
\begin{pro} \label{s-contraction}
    Let $(M,g)$ be a compact surface without conjugate points of genus greater than one, $\theta \in T_{1}M$ and $\mc{P}_N: D_N\subset S^s(U^s(\theta),\delta)\longrightarrow S^s(U^s(\theta),\delta)$ be as Proposition \ref{contraction}. Then, there exists $N'\geq 1$ such that for every $n\geq N'$, 
    \[ H_{u}\circ\mathcal{P}_n: V_{\delta}(\theta)\longrightarrow  V_{\delta}(\theta) \]
    is a topological contraction, i.e., there exists $\delta'<\delta$ such that $H_{u}\circ \mathcal{P}_n(V_{\delta}(\theta)) \subset V_{\delta'}(\theta)$.
\end{pro}

\begin{proof}
Let us recall the relative open neighborhood of $\mc{I}(\theta)$ in $S^{s}(U^s(\theta),\delta)$,
\[ W(\theta,\delta) = \bigcup_{\eta\in V_{\delta}(\theta)} \mc{W}^u(\eta). \]
Denote by $K$ the closure of $\mc{F}^s(\theta)\cap W(\theta,\delta)$ and by $U$ a relative open neighborhood of $\mc{I}(\theta)$ in $\mc{F}^s(\theta)$ such that $U\subset K$ strictly. By Lemma \ref{crossed-int}(2) we know that $K$ is a compact subset of $\mc{F}^s(\theta)$ containing a relative open neighborhood $U$ of $\mathcal{I}(\theta)$. Lemma \ref{contraction} provides $N'>0$ such that for every $n \geq N'$, 
\begin{equation}\label{inc}
    \Pi^h\circ \mc{P}_n(K)\subset U\subset K.
\end{equation}
So, for every $n\geq N'$, we can write the following expressions
\[ W_n=\bigcup_{\xi\in \Pi^h\circ \mc{P}_n(K)}\mc{W}^u(\xi), \quad  W=\bigcup_{\xi\in U}\mc{W}^u(\xi) \quad \text{ and } \quad  \overline{W}(\theta,\delta)=\bigcup_{\xi\in K}\mc{W}^u(\xi). \]
Thus for every $n\geq N'$, Equation \eqref{inc} implies that $W_n\subset W \subset \overline{W}(\theta,\delta)$ and hence 
\begin{equation}\label{inc2}
    W_n\cap V_{\delta}(\theta)\subset W\cap V_{\delta}(\theta)\subset \overline{W}(\theta,\delta)\cap V_{\delta}(\theta)=\overline{V}_{\delta}(\theta),
\end{equation}
where the second inclusion is strict. We also observe that for every $n\geq N'$, $W_n\cap V_{\delta}(\theta)=H_u\circ \mc{P}_n(V_{\delta}(\theta))$. Thus, Equation \eqref{inc2} ensures the existence of $\delta'<\delta$ such that for every $n\geq N'$,
$H_u\circ \mc{P}_n(V_{\delta}(\theta))\subset V_{\delta'}(\delta)\subset V_{\delta}(\theta)$ which concludes the proof.

\end{proof}

\section{Density of expansive points}

The goal of the section is to show Theorem \ref{main}, namely, expansive points are dense in $T_{1}M$ for every compact surface $(M,g)$ without conjugate points of genus greater than one. We follow the notations of previous sections. Let us begin by the following statement:
\begin{pro} \label{no-I-fol}
Let $(M,g)$ be a compact surface without conjugate points of genus greater than one, $\theta \in T_{1}M$, $U^s(\theta)$, $\delta>0$ and $S^s(U^s(\theta),\delta)$ as in Lemma \ref{Fol-neigh}. Then, there is no relative open neighborhood $\Sigma$ of $\theta$ in the cross section $S^s(U^s(\theta),\delta)$, composed only by pieces of non-trivial classes:
\begin{equation}\label{open}
    \Sigma=\bigcup_{\eta}J(\eta), \quad J(\eta)\subset \mc{I}(\eta), \quad \mc{I}(\eta) \text{ non-trivial,}
\end{equation}
where $J(\eta)$ is the maximal subset included in $\Sigma\cap \mc{I}(\eta)$.
\end{pro}

\begin{proof}
By contradiction, suppose that such a set $\Sigma$ exists for some $\theta \in T_{1}M$. From this we see that
\[ B= \bigcup_{|t| < \tau} \phi_{t}(\Sigma)\]
is an open set in $T_{1}M$ containing $\theta$. By Theorem \ref{v1}(4), there exists a periodic point $\theta'\in B$ very close to $\theta$. Since we can take a cross section around $\theta'$, let us suppose without loss of generality that $\theta$ is a periodic point with minimal period $T_{\theta}$. So, we can take the sequence of times $t_n=nT_{\theta}$ with $n\geq 1$ to get the sequence of Poincaré-type maps 
\[  \mc{P}_{nT_{\theta}}: D_{nT_{\theta}}\subset S^s(U^s(\theta),\delta)\longrightarrow S^s(U^s(\theta),\delta). \]
defined in Section \ref{scont}. The periodicity of $\theta$ implies that $\mc{P}_{nT_{\theta}}=\mc{P}_{T_{\theta}}^n$ for every $n\geq 1$. Since $\Sigma\subset S^s(U^s(\theta),\delta)$ is a relative open set in $S^s(U^s(\theta),\delta)$, there must exists a periodic point $\eta\in \Sigma$ such that $H_u(\eta)\in V_{\delta}(\theta)$. According to Lemma \ref{s-contraction}, there exists $N\geq 1$ such that $H_u\circ\mc{P}^N_{T_{\theta}}$ is a topological contraction. Thus, for $k\geq 1$ we get a sequence $H_u\circ\mc{P}^{kN}_{T_{\theta}}(\eta)$ converging to $\theta$ as $k\to \infty$. In particular, we deduce that $\mc{P}^{kN}_{T_{\theta}}(\eta)\in \Sigma$ for every $k\geq 1$. Moreover, there exists a subsequence $\mc{P}^{k_mN}_{T_{\theta}}(\eta)$ converging to some $\xi\in \mc{I}(\theta)$ as $m\to \infty$. Denote by $O(\eta)$ the orbit of $\eta$ by the geodesic flow. We note that for every $m\geq1$, $\mc{P}^{k_mN}_{T_{\theta}}(\eta)\in O(\eta)$. Since $\eta$ is periodic, $O(\eta)$ is compact and hence $\xi\in O(\eta)$. This is a contradiction because $\eta$ and $\xi$ have different orbits. 
\end{proof}

\begin{proof}[\bf Proof of Theorem \ref{main}]
By contradiction, suppose that expansive set is not dense in $T_1M$ hence there exists an open set $\Gamma$ composed only of non-expansive points. Choose $\theta\in \Gamma$, $\delta,\tau>0$ and $U^s(\theta)$ as in Lemma \ref{Fol-neigh} such that $\Gamma\subset\Gamma^s(U^s(\tilde{\theta}),\delta,\tau)$ (reducing $\Gamma$ if necessary). It is clear that $\Sigma=\Gamma\cap S^s(U^s(\theta),\delta)$ is a relative open neighborhood of $\theta$ in $S^s(U^s(\theta),\delta)$ composed only of non-trivial classes. Proposition \ref{tops}(4) implies that for each connected component of $\Sigma\cap \mathcal{F}^s(\eta)$ is composed by a subset of a single non-trivial class (for each component). Thus, $\Sigma$ has the form of Equation \eqref{open}, a contradiction.       
\end{proof}
 
\section{Some consequences of density of expansive dynamics}
This section is devoted to some applications of Theorem \ref{main} to some dynamical properties of the geodesic flow. Let us recall some objects related to the geodesic flow.

We outline the construction of the quotient flow of Section 4 of \cite{gelf19} in our context. For every $\eta,\theta\in T_1M$, $\eta$ and $\theta$ are equivalent, 
\[  \eta\sim\theta \quad \text{ if and only if }\quad  \eta\in \mathcal{I}(\theta).\]
This is an equivalence relation that induces a quotient space $X$ and a quotient map $\chi:T_1M\to X$. For every $\theta \in T_1M$, we denote by $[\theta]=\chi(\theta)$ the equivalence class of $\theta$. Using the geodesic flow $\phi_t$ induced by $(M,g)$, we define a quotient flow $\psi_t:X\to X$ by $\psi_t[\theta]=[\phi_t(\theta)]$ for every $t\in \mathbb{R}$. We shall endow $X$ with the quotient topology. A subset $A\subset T_1M$ is saturated with respect to $\chi$ if $A=\chi^{-1}\circ\chi(A)$. 

We next state a consequence of Lemma 3.1 of \cite{mam23} that has to do with a way of constructing open sets in $T_1M$ whose quotients are also open. 
\begin{lem}\label{sat}
Let $(M,g)$ be a compact manifold without conjugate points and with visibility universal covering $\tilde{M}$ and $\chi:T_1M\to X$ be the quotient map. Then, for every open set $U$ in $T_1M$ there exists an open saturated set $U'$ contained in $U$. In particular, $\chi(U')$ is an open set in $X$.
\end{lem}
We now give an application of Theorem \ref{main} concerning the full support of the maximal measure of the geodesic flow. First, let us give some properties of the quotient flow given in Section 6 and 7 of \cite{mam23}.
\begin{lem}\label{dyn}
Let $(M,g)$ be a compact surface without conjugate points of genus greater than one, $\mu$ be the maximal measure of the geodesic flow, $\chi:T_1M\to X$ be the quotient map, $\psi_t$ be the quotient flow. Then, $\psi_t$ has the specification property and has a unique maximal measure $\nu$ which has full support satisfying $\chi_*\mu=\nu$.  
\end{lem}
\begin{cor}\label{den}
If $(M,g)$ is a compact surface without conjugate points of genus greater than one then its maximal measure $\mu$ has full support.
\end{cor}
\begin{proof}
Let $U$ be an open set of $T_1M$. By Lemma \ref{sat}, there exists an open saturated set $U'\subset U$ in $T_1M$ such that $\chi(U')$ is open in $X$. We highlight that $U'$ is empty whenever $U$ does not contain any expansive point. However, Theorem \ref{main} provides an expansive point $\xi\in U$ which ensures that $U'$ is non-empty. By Lemma \ref{dyn}, the unique maximal measure $\nu$ has full support and hence $\nu(\chi(U'))>0$. The desired conclusion follows from
\[  \mu(U)\geq \mu(U')=\mu(\chi^{-1}\circ\chi(U'))=\chi_*\mu(\chi(U'))=\nu(\chi(U'))>0.  \]
\end{proof}
It is already known that the maximal measure of the geodesic flow has full support by a construction of the maximal measure done by Climenhaga-Knieper-War using the Patterson-Sullivan measure \cite{clim21}. However, the later point of view gives no information about neither the abundance of expansive points in $T_1M$ nor their relationship with the support of the maximal measure. Corollary \ref{den} gives the full support of the maximal measure without using Patterson-Sullivan theory. 

We turn to another application of Theorem \ref{main} that has to do with an alternative proof of the density of periodic points using the quotient flow.
\begin{cor}
Let $(M,g)$ be a compact surface without conjugate points of genus greater than one, $\chi:T_1M\to X$ be the quotient map and $\psi_t$ be the quotient flow. Then, periodic points are dense in $T_1M$.
\end{cor}
\begin{proof}
Let $U$ be an open set in $T_1M$. By Lemma \ref{sat}, there exists an open saturated set $U'\subset U$ such that $\chi(U')$ is open in $X$. As in proof of Corollary \ref{den}, Theorem \ref{main} provides an expansive point $\xi\in U$ which ensures that $U'$ is non-empty. Choose some piece of $\psi_t$-orbit $J$ intersecting $\chi(U')$. By the specification property of Lemma \ref{dyn}, there is a periodic orbit of $\psi_t$ shadowing $J$ and hence intersecting $\chi(U')$ at $\chi(\xi)$. Since $U'$ is saturated, we see that $\chi^{-1}(\chi(\xi))=\mc{I}(\xi)\subset U'$. If $T$ is the period of $\chi(\xi)$ then $\phi_T(\mc{I}(\xi))=\mc{I}(\xi)$. Thus, Brouwer's fixed point Theorem implies that $\mc{I}(\xi)$ has a fixed point $\eta$ which is a periodic point of the geodesic flow. This concludes the proof since $\eta\in U'\subset U$. 
\end{proof}

\section{Expansiveness of the geodesic flow versus expansive points}
The purpose of the section is to discuss the relationship between the usual notion of expansivity for a flow and the notion of expansive point in the setting of geodesic flows. 

It is quite natural to expect that topological transversality of the submanifolds $\mathcal{F}^{s}(\theta)$ and $\mathcal{F}^{u}(\theta)$ for every $\theta\in T_{1}M$ characterizes the expansiveness of the geodesic flow. There is a similar characterization for Anosov geodesic flows. In the context of compact manifolds without conjugate points, we would like to show that expansiveness of the geodesic flow is actually almost equivalent to the expansivity of every point in the unit tangent bundle. 

Recall that for a compact manifold without conjugate points, the expansive set of its geodesic flow is defined by 
\[  \mc{R}_0=\{ \xi\in T_1M: \mc{F}^s(\xi)\cap \mc{F}^u(\xi)=\{ \xi\} \}. \]
\begin{pro}\label{exp}
Let $(M,g)$ be a compact manifold without conjugate points and $\T{M}$ be its universal covering. If $\T{M}$ is quasi-convex where geodesic rays diverge then the geodesic flow restricted to $\mc{R}_0$ is expansive.   
\end{pro}
\begin{proof}
Let $\epsilon>0$ be smaller than half of the injectivity radius of $T_{1}M$ endowed with Sasaki metric. Let $\xi,\xi'\in \mc{R}_0$ and $s:\mathbb{R}\to \mathbb{R}$ be a time reparametrization such that $d_s(\phi_t(\xi),\phi_{s(t)}(\xi'))\leq \epsilon$ for every $t\in \mathbb{R}$. Recall that $\pi: \tilde{M} \to M$ is the covering map and $d\pi : T_{1}\tilde{M} \to T_{1}M$ is the induced covering map.

By choice of $\epsilon$ and elementary topology arguments, for each lift $\T{\xi}\in T_1\T{M}$ of $\xi$ there exists a lift  $\T{\xi}'\in T_1\T{M}$ of $\xi'$ such that 
\[ d_s(\T{\phi}_t(\T{\xi}),\T{\phi}_{s(t)}(\T{\xi}'))\leq \epsilon, \quad \text{ for every }t\in \mathbb{R}.  \]  
Thus, the Hausdorff distance between the orbits of $\T{\xi}$ and $\T{\xi}'$ is bounded above. This together with Theorem \ref{asim} yields that $\T{\xi}'\in \mc{\T{I}}(\T{\phi}_u(\T{\xi}))$ for some $u\in \mathbb{R}$. Applying the covering map $d\pi$ we see that $\xi'\in \mc{I}(\phi_u(\xi))=\phi_u(\xi)$ because $\xi$ is an expansive point. We can choose $\epsilon$ small enough so that $|u|< \epsilon$.
\end{proof}
Let us see some consequences of Proposition \ref{exp}. The next statement involves the category of manifolds without conjugate points whose Green bundles are continuous. Since this is the only place of the article where Green bundles appear, we shall refer the reader to \cite{gelf20} for the definition and further properties. 
\begin{pro}\label{mea}
Let $(M,g)$ be a compact manifold without conjugate points, $\phi_t$ be its geodesic flow and $\T{M}$ be its universal covering. Assume that $\T{M}$ is quasi-convex and geodesic rays diverge.
\begin{enumerate}
    \item If $\T{M}$ is a visibility manifold and $M$ has a hyperbolic closed geodesic and continuous Green bundles then $\phi_t$ is expansive on a dense open set of $T_1M$.
    \item If $M$ is a surface then $\phi_t$ is expansive on the dense set $\mc{R}_0$.
    \item If $\mu$ is a $\phi_t$-invariant measure giving zero measure to non-expansive points then $\phi_t$ is expansive $\mu$-almost everywhere.
    \item If $\mc{R}_0=T_1M$ then $\phi_t$ is expansive.  
\end{enumerate}
\end{pro}
Note that $\mc{R}_0=T_1M$ means topological transversality of the horospheres $\mc{F}^s(\theta)$ and $\mc{F}^u(\theta)$ at $\theta$ for every $\theta\in T_1M$.
\begin{proof}
For item 1, from Theorem 1.3(3) of \cite{mam24} it follows that $\mathcal{F}^{s}(\theta)$ and $\mathcal{F}^u(\theta)$ are always tangent to Green bundles. Moreover, Green bundles are continuous and transverse in an open neighborhood of a hyperbolic closed geodesic. In the particular case of item 1, the geodesic flow is topologically transitive. As in Corollary 8.1 of of \cite{mam24}, it follows that $\mc{R}_1=\{\xi\in T_1M:G^s(\xi)\cap G^u(\xi)=\{\xi\} \}$ is an open dense set contained in $\mc{R}_0$ where $G^s(\xi)$ and $G^u(\xi)$ are the Green subspaces at $\xi$. Thus, Proposition \ref{exp} proves item 1. Theorem \ref{main} and Proposition \ref{exp} show item 2. Item 3 follows from Proposition \ref{exp} noting that support of $\mu$ is included in $\mc{R}_0$. Item 4 follows straightforward from Proposition \ref{exp}.
\end{proof}
We finish the article with a result that almost provides the equivalence between expansivity of a geodesic flow and expansivity of every point in the unit tangent bundle. 
\begin{teo} \label{Exp0}
Let $(M,g)$ be a compact manifold without conjugate points and $\tilde{M}$ be its universal covering. Then, the geodesic flow is expansive if and only if $\tilde{M}$ is quasi-convex, geodesic rays diverge and $\mc{R}_0=T_1M$. 
\end{teo}
For the proof we need the following alternative characterization of expansive geodesic flows \cite{rugg94,rugg97}. 
\begin{teo} \label{Exp1}
Let $(M,g)$ be a compact manifold without conjugate points and $\tilde{M}$ be its universal covering. Then, the geodesic flow is expansive if and only if:
	\begin{enumerate}
		\item $\tilde{M}$ is a visibility manifold. In particular, $\tilde{M}$ is quasi-convex and geodesic rays diverge uniformly. 
		\item For every $\theta \in T_{1}\tilde{M}$, $\mathcal{\T{F}}^{s}(\theta)$ and $\mathcal{\T{F}}^{u}(\theta)$ are the strong stable and unstable sets of $\theta$. Namely, for every $D,\epsilon>0$ there exists $T(D,\epsilon)>0$ such that for every $\eta\in\mathcal{\T{F}}^{s}(\theta)$ with $d_s(\theta,\eta)\leq D$,
		\[ d_s(\T{\phi}_{t}(\theta),\T{\phi}_{t}(\eta))\leq \epsilon, \quad \text{for every } t \geq T(D,\epsilon) \] 
		  and for every $\eta\in\mathcal{\T{F}}^{u}(\theta)$ with $d_s(\theta,\eta)\leq D$,
		\[ d_s(\T{\phi}_{t}(\theta),\T{\phi}_{t}(\eta))\leq \epsilon, \quad \text{for every } t \leq -T(D,\epsilon). \]  
	\end{enumerate}
\end{teo}

\begin{proof}[Proof of Theorem \ref{Exp0}]
The reverse implication is just item 4 of Proposition \ref{mea}. For the direct implication, by contradiction suppose the geodesic flow is expansive but $\mc{R}_0\neq T_1M$. Hence there exist a non-expansive point $\xi\in T_1M$ and a non-trivial class $\mc{I}(\xi)\subset T_1M$. Let $\tilde{\xi}\in T_1\T{M}$ be any lift of $\xi$. By Theorem \ref{Exp1}(1) and Theorem \ref{connected-class}, we know that 
\[  \mathcal{\T{I}}(\tilde{\xi})= \mathcal{\T{F}}^{s}(\tilde{\xi}) \cap \mathcal{\T{F}}^{u}(\tilde{\xi}) \]
is a compact connected set with diameter bounded by $C>0$. Now, given $\epsilon>0$ let $T(C, \epsilon)>0$ be as in Theorem \ref{Exp1}(2). Then, for every $\tilde{\eta} \in \mathcal{\T{I}}(\tilde{\xi})$ we have  
	\[   d_s(\T{\phi}_{t}(\tilde{\eta}), \T{\phi}_{t}(\tilde{\xi})) \leq \epsilon, \quad \text{ for } |t|\geq T(C, \epsilon). \]
By continuity of the geodesic flow upon initial conditions, there exists $\delta>0$ such that if $\tilde{\eta}\in \mathcal{\T{I}}(\tilde{\xi})$ with $d_s(\tilde{\eta}, \tilde{\xi}) \leq \delta$ then 
\[  d_s(\T{\phi}_{t}(\tilde{\eta}), \T{\phi}_{t}(\tilde{\xi})) \leq \epsilon, \quad \text{ for }|t|\leq T(C,\epsilon). \] 
Since $\mathcal{\T{I}}(\tilde{\xi})$ is connected, it meets every suitably small ball around $\tilde{\xi}$. Therefore, 
\[  d_s(\T{\phi}_{t}(\tilde{\xi}), \T{\phi}_{t}(\tilde{\eta})) \leq \epsilon, \quad \text{ for every } t\in \mathbb{R}. \] 
Since $\epsilon>0$ was arbitrary, the geodesic flow cannot be $\epsilon$-expansive for any $\epsilon >0$, which is a contradiction. 
\end{proof}
Since compact higher genus surfaces without conjugate points have quasi-convex universal coverings and geodesic rays diverge uniformly, we obtain a simplified characterization of expansive geodesic flows from Theorem \ref{Exp0}.
\begin{cor}
Let $(M,g)$ be a compact surface without conjugate points of genus greater than one. Then, the geodesic flow is expansive if and only if $\mc{R}_0=T_1M$.
\end{cor}
Thus, in this context expansivity is completely characterized by topological transversality in all stable and unstable horospheres.
\bigskip

\textbf{Conjecture:} The reverse implication in Theorem \ref{Exp0} holds without the assumptions of quasi-convexity and divergence of geodesic rays. 

	\bibliographystyle{plain}
	\bibliography{ref}
	
\end{document}